\documentclass[12pt]{amsart}
\usepackage{amssymb}
\usepackage{amsmath}
\usepackage{verbatim}
\usepackage{pictex}

\usepackage              {hyperref} 
\usepackage[initials,nobysame,author-year, 
 non-sorted-cites]       {amsrefs}  
\usepackage[english]     {babel}    
\usepackage        [utf8]{inputenc} 
\usepackage          [T1]{fontenc}  

\theoremstyle{plain}
 \newtheorem {theorem}{Theorem}[section]
 \newtheorem {proposition}{Proposition}[section]
 \newtheorem {lemma}{Lemma}[section]
\theoremstyle{definition}

\makeatletter \catcode`\'=11 
\let\print@citenames\CiteNamesFull 
\renewcommand\@biblist[1][]{
    \stepcounter{bib@env}
    \biblistfont
    \labelsep 2em\relax 
    \let\@bibitem\amsrefs@bibitem
    \let\@lbibitem\amsrefs@lbibitem
    \list{\BibLabel}{%
        \restore@labelwidth
        \itemindent=-\leftmargin 
        \@maxlabelwidth\z@
        \@nmbrlisttrue
        \def\@listctr{bib}%
        \let\makelabel\bib@mklab
        #1\relax
    }%
    \sloppy
    \interlinepenalty\@m
    \clubpenalty\@M
    \widowpenalty\clubpenalty
    \frenchspacing
    \ResetCapSFCodes
    \@ifstar{\@biblistsetup}{}%
}
\def\strip@Zblprefix#1#2#3#4#5\@nil{%
    \def\@tempa{#1#2#3#4#5}%
    \if#1Z%
        \if#2b%
        \if#3l%
            \def\@tempa{#4#5}%
        \fi
        \fi
    \fi
}
\def\Zbl#1{%
    \relax\ifhmode\unskip\spacefactor3000 \space\fi
    \begingroup
        \strip@Zblprefix#1\@nil
        \edef\@tempa{Zbl~\@nx\Zblhref{\@tempa}{\@tempa}}%
    \@xp\endgroup
    \@tempa
}
\providecommand{\Zblhref}[2]{#1}
\def\strip@JFMprefix#1#2#3#4#5\@nil{%
    \def\@tempa{#1#2#3#4#5}%
    \if#1J%
        \if#2F%
        \if#3M%
            \def\@tempa{#4#5}%
        \fi
        \fi
    \fi
}
\def\JFM#1{%
    \relax\ifhmode\unskip\spacefactor3000 \space\fi
    \begingroup
        \strip@JFMprefix#1\@nil
        \edef\@tempa{JFM~\@nx\JFMhref{\@tempa}{\@tempa}}%
    \@xp\endgroup
    \@tempa
}
\providecommand{\JFMhref}[2]{#1}
\makeatother \catcode`\'=12 

\newcounter{firstpagenumber}
\begin{document}

\title
{Diffusion on homogeneous ultrametric spaces:\\ the contributions of 
Alessandro Figà-Talamanca}

\author{\bf Wolfgang Woess}
\address{\parbox{.8\linewidth}{Institut f\"ur Diskrete Mathematik,\\ 
Technische Universit\"at Graz,\\
Steyrergasse 30, A-8010 Graz, Austria\\}}
\email{\rm woess@tugraz.at}
\date{January 6, 2026} 
\begin{abstract}
Alessandro Figà-Talamanca (1938-2023) was an influential Italian mathematician,
scientific leader of the Italian group of harmonic analysis for many years.
Since the late 1970ies, his interest focussed on harmonic
analysis on free groups and trees. In the later years of his scientific
work he became also interested in diffusion processes on homogeneous 
ultrametric spaces such as local fields and totally disconnected Abelian groups.
This is related with the close connection of those spaces with trees and their 
boundaries and concerns, in particular, the construction of such processes 
via discrete-time  walks on trees. 
The present notes provide rather detailed comments on this
part of his work and the related, quite abundant literature. 
This is intended to become part of a volume of selected papers by
Figà-Talamanca, accompanied by comments such as the present text.
\end{abstract}


\maketitle

\section{Introduction}\label{sec:intro}

\medskip 

These notes provide comments on the following 
articles of A. Figà-Talamanca and coauthors plus related literature.
%
\smallskip

{
}

\bigskip

The late work of Alessandro Figà-Talamanca concerns  
diffusion processes on ultrametric spaces and their description via  
random walks on trees. It started with \ocite{FT-41}, \ycite{FT-44} for compact spaces.
The main contents of these two papers coincide. Here, we refer mainly to the 
first one.\footnote{The second one was written shortly before the first one,
as the notes of a minicourse at a conference, but its proceedings volume 
appeared much later.} The first part of \ocite{FT-43} applies similar methods  
with a slightly different goal in a non-compact setting.
The two main papers 
are \ocite{FT-46}, \ycite{FT-47}.
There is an ample and diverse literature on the subject of diffusion and 
Laplace-type operators on ultrametric spaces, with considerably varying notation
and many mutual citations which are not necessarily based on mutual reading.
As the author of the present notes who was only relatively briefly involved in this 
fascinating field of work via 
\ocite{Bendikov-et-al2014} and \ocite{Bendikov-et-al2019}, I will try here to 
elaborate a common picture, connecting the mentioned papers of Figà-Talamanca
with a selection from the large body of other important work and translating different
notation and terminolgy into one setting. Since it is not expected
that every potential reader is acquainted with the relation between ultrametric spaces 
and trees, the exposition is quite detailed, starting with the basics.


\section{Ultrametric spaces and trees}\label{sec:trees}


An \emph{ultrametric space} is a metric space $(\mathcal{X},\rho)$ which in place of the triangle equality
satisfies the stronger \emph{utrametric inequality}
$$
\rho(x,z) \le \max \{ \rho(x,w), \rho(w,z) \}.
$$
As a consequence, if $y \in B(x,r)$, the closed ball with centre $x$ and radius $r$,
then $B(x,r) = B(y,r)$. Closed balls are either disjoint or one contains the other.
We always assume that the space is separable and that closed balls are compact.
This implies that for any $x \in \mathcal{X}$, the set 
$$
\Lambda(x) = \{ \rho(x,y) : y \ne x \}
$$ 
does not accumulate in $(0\,,\,\infty)$, whence  
$
\Lambda = \bigcup_{x \in \mathcal{X}} \Lambda(x)
$
is countable. Here we also assume that $\mathcal{X}$ has no isolated points. 

The classical example is the field of $\mathsf{p}$-adic numbers, where $\mathsf{p}$ is a prime.
Before further examples, we explain the well-known realisation of ultrametric spaces
as boundaries of trees. There are many references for this basic equivalence, e.g. 
\ocite{Rammal}, also \ocite{FT-42}.

Let $T$ be an infinite, locally finite tree with set of non-oriented edges 
$E = E_T$.  We write $u, v,\dots$ for its vertices and $[u,v] \equiv [v,u]$ for edges, in which we also 
write $u \sim v$, the neighbourhood relation.  The \emph{degree} $\deg(u)$ of a vertex is
the number of its neighbours. The graph distance $d(u,v)$ is the length (number of edges)
of the unique minimal path (geodesic arc) $\pi(u,v)$ between the vertices $u$ and $v$.
A \emph{geodesic ray} in $T$ is a one-sided infinite path 
$[u_0\,,u_1\,,u_2\,,\dots]$ where $u_i \sim u_{i-1}$ and no vertex appears twice.
Two geodesic rays are called equivalent, if their symmetric difference is finite, i.e.,
they differ by finite initial pieces. An \emph{end} of $T$ is an equivalence class of rays,
and the set of all ends of $T$ is the \emph{boundary} $\partial T$, later to be 
denoted $\mathcal{X}$. We write $x,y,z$ etc for the elements of $\partial T$. 
For every $u \in T$ and $x \in \partial T$ there is precisely one geodesic ray $\pi(u,x)$
which starts at $u$ and represents $x$. 
\\[5pt]
\textbf{A) Compact case.} 
We choose a reference vertex (origin) $o \in T$. For any vertex $u$ we write $|u|$ for its 
graph distance from $o$. If $u  \ne o$ then the \emph{predecessor} $u'$ of $u$ is 
the neighbour of $u$ on $\pi(o,u)$. If $x,y \in T \cup \partial T$ are distinct, then 
$x \wedge y$ denotes the
common vertex on $\pi(o,x) \cap \pi(o,y)$ which is furthest from $o$, that is,
$$
\pi(o,x) \cap \pi(o,y) = \pi(o, x \wedge y).
$$
We now choose any real number $q > 1$ and define
$$
\rho(x,y) = \begin{cases} 0\,,&\text{if}\;x=y\\
                       q^{-|x \wedge y|}\,,&\text{if}\;x\ne y.   
         \end{cases}
$$
This makes $\mathcal{X} = \partial T$ into a compact, separable ultrametric
space. Closed (open-compact) balls in this metric are as follows: let $u \in T$, and 
let 
$$
\partial T_u = \{ x \in \partial T : \pi(o,x) \;\text{passes through}\; u \}.
$$
This is the closed ball $B(x,q^{-n})$, where $n = |x|$. In particular, the induced 
topology does not depend on the choice of $q$, and there are various other ways to define
a topologically equivalent ultrametric with the same collection of compact-open balls. 
We also see the following: if there is $u \in T \setminus \{o\}$ with $\deg(u)=2$ then
$\partial T_u = \partial T_{u'}$, so that this ball of radius $q^{-|u|}$ also has
radius $q^{-|u|+1}$. We avoid this: for the compact case, we assume that $\deg(u) \ge 3$
for all $u \in T \setminus \{o\}$. 

\smallskip

Conversely, start with a compact ultrametric space $\mathcal{X}$ without isolated points. 
Consider the set 
$$
T = \{ u = B(x,r) : x \in \mathcal{X}, r \in \Lambda(x)\}.
$$
It is countable. Let $r_0 = \max \Lambda = \max \Lambda(x)$, so that 
$o := B(x,r_0) = \mathcal{X}$ for every $x$. Recall that balls are open-compact. 
Each ball $u = B(x,r) \in T \setminus \{o\}$  has a unique predecessor 
ball $u' = B(x,r')$ where $r' = \min \Lambda(x) \cap (r\,,\,\infty)$, and 
the number of successors of $u$, i.e., the balls whose predecessor is $u$, is finite. 
It is always  $\ge 2$, since we assume that $\mathcal{X}$ has no isolated points.
This gives $T$ the structure of an infinite tree where $\deg(u) \ge 3$ for all 
vertices $u \ne o$.
Any end of $T$ is described by an infinite, decreasing 
sequence of balls. Its intersection is a unique $x \in \mathcal{X}$ which is a 
centre of each of those balls. If $x, y \in \mathcal{X}$ are distinct, then in 
$T$ their confluent $u=x \wedge y$ is $u = B(x,r) = B(y,r)$, where $r = \rho(x,y)$ 
(not necessarily of the form $q^{-|x \wedge y|}$, according to the choice of the 
ultrametric).

For later use, we select one elemement of $\mathcal{X} \equiv \partial T$
which we denote  $\mathit{0}$ (later on to become a neutral element) and write 
$\pi(o,\mathit{0}) = [o=o_0\,,o_1\,,o_2\,,\dots]$.

\begin{figure}[b!]
\beginpicture 

\setcoordinatesystem units <.42mm,.70mm> 

\setplotarea x from -250 to 106, y from -20 to 80


\arrow <6pt> [.2,.67] from -198 2 to -205 -5

\plot -198 2 -136 64 /

\plot -168 32 -138 2 /

 \plot -184 16 -170 2 /

 \plot -152 16 -166 2 /

 \plot -192 8 -186 2 /

 \plot -176 8 -182 2 /

 \plot -160 8 -154 2 /

 \plot -144 8 -150 2 /

 \plot -196 4 -194 2 /

 \plot -188 4 -190 2 /

 \plot -180 4 -178 2 /

 \plot -172 4 -174 2 /

 \plot -164 4 -162 2 /

 \plot -156 4 -158 2 / 

 \plot -148 4 -146 2 /

 \plot -140 4 -142 2 /

 \plot -74 2  -136 64 /

 \plot -134 2 -104 32 /

 \plot -130 2 -132 4 /

 \plot -126 2 -124 4 /

 \plot -122 2 -128 8 /

 \plot -118 2 -112 8 /

 \plot -114 2 -116 4 / 

 \plot -110 2 -108 4 /

 \plot -106 2 -120 16 /
 
 \plot -102 2 -88 16 /
 
 \plot -98 2 -100 4 /
 
 \plot -94 2 -92 4 /
 
 \plot -90 2 -96 8 /
 
 \plot -86 2 -80 8 /

 \plot -82 2 -84 4 /
 
 \plot -78 2 -76 4 /

\put {$\vdots$} at -168 -3
\put {$\vdots$} at -104 -3

\put {$S_1$} [l] at -226 32
\put {$S_2$} [l] at -226 16
\put {$S_3$} [l] at -226 8.7
\put {$S_4$} [l] at -226 4
\put {$\partial T$} [l] at -227 -12 
\put {$\vdots$} at -220 -3
\put {$\mathit{0}$} at -208 -8

\put {$o\!=\!o_0$} at -135 66
\put {$o_1$} [r] at -168 35 
\put {$o_2$} [r] at -184 19
\put {$u$} [r] at -120 19
\put {$v'$} [l] at -88 19.5
\put {$v$} [l] at -80 10.5
\put {$u \!\wedge\! v$} [l] at -104.5 35

\put {$\scriptstyle \bullet$} at -104 32
\put {$\scriptstyle \bullet$} at -184 16
\put {$\scriptstyle \bullet$} at -168 32
\put {$\scriptstyle \bullet$} at -80 8
\put {$\scriptstyle \bullet$} at -120 16
\put {$\scriptstyle \bullet$} at -88 16
\put {$\scriptstyle \bullet$} at -136 64


\setlinear 

\arrow <6pt> [.2,.67] from 2 2 to 80 80

\arrow <6pt> [.2,.67] from 2 2 to -5 -5

\plot 32 32 62 2 /

 \plot 16 16 30 2 /

 \plot 48 16 34 2 /

 \plot 8 8 14 2 /

 \plot 24 8 18 2 /

 \plot 40 8 46 2 /

 \plot 56 8 50 2 /

 \plot 4 4 6 2 /

 \plot 12 4 10 2 /

 \plot 20 4 22 2 /

 \plot 28 4 26 2 /

 \plot 36 4 38 2 /

 \plot 44 4 42 2 /

 \plot 52 4 54 2 /

 \plot 60 4 58 2 /

 \plot 99 29 64 64 /

 \plot 66 2 96 32 /

 \plot 70 2 68 4 /

 \plot 74 2 76 4 /

 \plot 78 2 72 8 /

 \plot 82 2 88 8 /

 \plot 86 2 84 4 /

 \plot 90 2 92 4 /

 \plot 94 2 80 16 /
 

\setdots <3pt>
\putrule from -213 4 to -76 4  
\putrule from -213 8 to -80 8  
\putrule from -213 16 to -88 16  
\putrule from -213 32 to -104 32  

\putrule from -10.8 4 to 102 4  
\putrule from -10.1 8 to 102 8  
\putrule from -8 16 to 102 16   
\putrule from -7.7 32 to 102 32   
\putrule from -7.7 64 to 102 64   

\linethickness .8pt
\setdashes <3pt>
\putrule from -212 -12 to -66 -12  
\putrule from -12 -12 to 109 -12   
\put {$\vdots$} at 32 -3
\put {$\vdots$} at 84 -3

\put {$\dots$} [l] at 103 6
\put {$\dots$} [l] at 103 48

\put {$H_{-2}$} [l] at -26 64
\put {$H_{-1}$} [l] at -26 32
\put {$H_{0}$} [l] at -26 16
\put {$H_1$} [l] at -26 8.7
\put {$H_2$} [l] at -26 4
\put {$\partial^*T$} [l] at -30 -12 
\put {$\vdots$} at -20 -3
\put {$\vdots$} [B] at -16 70
\put {$\varpi$} at 82 82
\put {$\mathit{0}$} at -8 -8

\put {$o_{-1}$} at 24 34
\put {$o\!=\!o_0$} at 4.5 18.5
\put {$o_1$} at 4 10.5
\put {$u$} at 57.8 10.5
\put {$u'$} at 53 19
\put {$v$} at 99 34.5
\put {$u \!\curlywedge\! v$} [r] at 63 66.5

\put {$\scriptstyle \bullet$} at 8 8
\put {$\scriptstyle \bullet$} at 16 16
\put {$\scriptstyle \bullet$} at 32 32
\put {$\scriptstyle \bullet$} at 56 8
\put {$\scriptstyle \bullet$} at 48 16
\put {$\scriptstyle \bullet$} at 96 32
\put {$\scriptstyle \bullet$} at 64 64

\endpicture
\end{figure}

\medskip

\noindent
\textbf{B) Non-compact case.} We start with a tree $T$ as above and assume that 
$\deg(u) \ge 3$ for each vertex. We choose one reference end $\varpi \in \partial T$ 
and consider the punctured boundary $\partial^*T = \partial T \setminus \{ \varpi \}$.
For each $x \in \partial^* T$, there is a unique bi-infinite geodesic
$\pi(\varpi,x) = [\dots, u_{-2}\,, u_{-1}\,,u_0\,,u_1\,,u_2\,,\dots]$
of vertices with $u_i \sim u_{i-1}$, 
such that 
$[u_k\,,u_{k+1}\,,u_{k+2}\,,\dots] = \pi(u_k\,,x)$ and  
$[u_k\,,u_{k-1}\,,u_{k-2}\,,\dots] = \pi(u_k\,,\varpi)$. 
We fix again one element $\mathit{0} \in \partial^*T$ and write $o_k\,$, 
$k \in \mathbb{Z}$, for 
the elements of $\pi(\varpi,\mathit{0})$, in particular $o = o_0\,$. 

For distinct $x,y \in T \cup \partial^* T$, we denote by $x \curlywedge y$ the first common
vertex on $\pi(x, \varpi)$ and $\pi(y,\varpi)$, that is,
$$
\pi(x\curlywedge y,\varpi) = \pi(x,\varpi) \cap \pi(y,\varpi).
$$
It is always well-defined. Also, $u \curlywedge u = u$ for
$u \in T$. For any vertex $u$, we define its \emph{horocycle number} 
$\mathsf{hor}(u) = d(u,u \curlywedge o) - d(o,u \curlywedge o)$. 
For $k\in \mathbb{Z}$, the \emph{$k^{\text{th}}$ horocycle} is 
$H_k = \{ u\in T : \mathsf{hor}(u)=k \}$.   
Every element $u$ of $H_k$ has precisely one neighbour in $H_{k-1}\,$ which, 
by abuse of the 
terminology in the compact case, we call the predecessor $u'$, and it has 
$\deg(u)-1$ successors
in $H_{k+1}\,$.
Again, we choose $q > 1$ and define the metric on $\partial^*T$ by 
\begin{equation}\label{eq:rho}
\rho(x,y) = \begin{cases} 0\,,&\text{if}\;x=y\\
                       q^{-\mathsf{hor}(x \curlywedge y)}\,,&\text{if}\;x\ne y.   
         \end{cases}
\end{equation}
This makes $\mathcal{X} = \partial^*T$ into a non-compact, separable ultrametric space
without isolated points in which balls are open-compact. Any ball is of the following form: 
for a vertex $u$, let $\partial^* T_u = \{ x \in \partial^* T : u \in \pi(\varpi,x) \}$.
Its radius $\equiv$ diameter with respect to the above metric 
$\rho$ is $q^{-\mathsf{hor}(u)}$. Again,  various other ultrametrics with the same 
collection of compact-open balls can be used.

\smallskip

Conversely, start with an ultrametric space $\mathcal{X}$ with the above properties.
Then we take again the collection of all closed (open-compact) balls. These form our
vertex set. By the properties of the set of distances $\Lambda$ and its subsets 
$\Lambda(x)$, each ball $B$ is contained in a unique ball $B'$ with the next bigger 
radius $\equiv$ diameter, its 
predecessor, and by compactness, it is the disjoint union of finitely many balls of which
it is the predecessor. This provides the tree structure.
Again, since there are no isolated points, each $x \in \mathcal{X}$ is the 
intersection of a descending sequence of balls, which provide a geodesic in the tree whose 
limit point (equivalence class of rays) is $x$.\footnote{We remark that isolated 
points can also be handled; they come up as leaves of the tree 
(vertices with degree $1$), to be included in the boundary.}

As a unified notation in the compact as well as the non-compact case, for a vertex 
$u\in T$ we write 
$$
\mathsf{q}(u) = \bigl| \{ v \in T : v' =u \}|
$$
for the forward degree, which is equal to $\deg(u)-1$ unless $u=o$ in the compact case,
where $\mathsf{q}(o)=\deg(o)$.

\smallskip 

In the basic non-compact example where $\mathcal{X} = \mathbb{Q}_{\mathsf{p}}\,$, 
the ring of $\mathsf{p}$-adic numbers (field when $\mathsf{p}$ is prime),
the associated tree is such that each vertex has exactly $\mathsf{p}$ successors. We choose
$\mathit{0}$ as the neutral element in $\mathbb{Q}_{\mathsf{p}}\,$. Then 
$\mathsf{hor}(x \wedge \mathit{0})$ is the 
valuation of $x \in \mathbb{Q}_{\mathsf{p}} \setminus \{\mathit{0}\}$. For the 
ultrametric of \eqref{eq:rho}, here one chooses the base $q = \mathsf{p}$.

\smallskip

The work of Figà-Talamanca addressed here concerns the situation when the space 
$\mathcal{X}$
is homogeneous, i.e., its isometry group acts transitively. For the associated tree, this
means that the group acts transitively on the level sets, which are
the spheres 
$S_k = \{ u \in T: |u|=k\}$   ($k \in \mathbb{N}$) in the compact case, resp. 
the horocycles  $H_k$  ($k \in \mathbb{Z}$) in the non-compact case. 
In particular, all vertices $u$ in the 
$k^{\text{th}}$ level sets have the same forward degree $\mathsf{q}(u)=\mathsf{q}_k\,$,
and as balls in $\mathcal{X}$, they have the same diameter $\equiv$ radius $r_k\,$.
Any transitive, closed 
subgroup of  the isometry group which acts freely on $\mathcal{X}$ is locally 
compact and totally disconnected,  and in the non-compact case, it is 
the union of a countable family of compact subgroups (in $T$,
these are the stabilisers of the vertices $o_{-k}$, $k \in \mathbb{N}$). 
Conversely, if $\mathcal{X}$ itself is such a group, then the associated 
tree has constant degrees on its level sets. In particular, there are 
\emph{Abelian} groups which realise any given sequence of forward 
degrees $\mathsf{q}_k \ge 2$. 
An important example is the non-compact group which is denoted 
$\Omega_{\mathbf{a}}$ in \ocite{Hewitt-Ross}. Here 
$\mathbf{a} = (\mathbf{a}_k)_{k \in \mathbb{Z}}$ is an arbitrary 
sequence of natural numbers $\ge 2$, and the associated tree is such that 
$\mathsf{q}_k = \mathbf{a}_k$ on $H_k$. It contains the compact subgroup denoted 
$\Delta_{\mathbf{a}}$, which in terms of the action on the tree is the stabiliser 
of the vertex $o= o_0$ in $\Omega_{\mathbf{a}}$, and corresponds to $\partial^* T_o$ 
in the boundary. When $\mathbf{a}_k = p \ge 2$ for all $k \in \mathbb{Z}$,
one has $\Omega_{\mathbf{a}} = \mathbb{Q}_{\mathsf{p}}\,$, the (additive group in the) 
ring of $\mathsf{p}$-adic numbers (field, when $\mathsf{p}$ is prime), see above, and 
$\Delta_{\mathbf{a}}$ is the group (ring) of $\mathsf{p}$-adic integers.


\section{Laplacians on ultrametric spaces}\label{sec:Lap}


The main purpose of 
\ocite{FT-41} and \ocite{FT-46}, \ycite{FT-47} is to construct
diffusion processes adapted to the structure of the space via random walks on 
the corresponding  tree whose boundary is identified with $\mathcal{X}$.  
The result is a semigroup $(P_t)_{t > 0}$ of transition operators with an
associated heat kernel $p_t(x,y)$ with respect to a natural
reference measure $\mathsf{m}$. To the semigroup, there corresponds its infinitesimal 
generator $-\mathfrak{L}$, where $\mathfrak{L}$ is a non-negative definite, typically
unbounded operator, often called the \emph{Laplacian} (or maybe minus the Laplacian) 
on $L^2(\mathcal{X}, \mathsf{m})$.
In all the relevant work, $\mathfrak{L}$ has a pure point spectrum with compactly 
supported  eigenfunctions. Thus, in order to see how the constructions of 
Del Muto and Figà-Talamanca are embedded in the larger body of work, we first 
provide an overview over the basic ultrametric Laplacians and a 
comparison of their spectral properties. Thereafter, in \S \ref{sec:RW}, 
we shall look at the approaches via random walks on trees.
\\[5pt]
\textbf{A) Hierarchical Laplacians.}
Before reviewing the developments in chronological order, we start with 
a more recent model which is elegant in view of its simplicity, the 
\emph{hierarchical Laplacians} $\mathcal L$, which appear as minus the infinitesimal
generators of \emph{isotropic Markov semigroups} on ultrametric spaces. 
It was developed systematically by the 
first two authors of 
\ocite{Bendikov-at-al2012} and \ocite{Bendikov-et-al2014}, see also the surveys by
\ocite{Bendikov2018} and \ocite{Grigoryan2023} and the references therein.
Without necessarily specifying the metric, we need the following.
\begin{itemize}
\item[(i)] The collection $\mathcal{B}$ of all open-compact balls 
(equivalently, the tree $T$), 
\item[(ii)] a non-atomic Radon measure $\mathsf{m}$ on $\mathcal{X}$ such that 
$\mathsf{m}(B) >0$ for each ball and $\mathsf{m}(\mathcal{X}) = \infty$ when 
$\mathcal{X}$ is non-compact,\footnote{In the non-compact case, the 
assumption $\mathsf{m}(\mathcal{X}) = \infty$ is sometimes dropped.} and
\item[(iii)] a \emph{choice function} $C: \mathcal{B} \to (0\,,\,\infty)$ such that for 
each ball $B$ and each $x \in \mathcal{X}$
$$
\lambda(B) = \sum_{\tilde B \in \mathcal{B}\,:\, B \subseteq \tilde B} C(\tilde B) < \infty
\quad \text{and}\quad \sum_{B \in \mathcal{B}\,:\, x \in B} C(B) = \infty \,.  
$$
\end{itemize}
In terms of $T$, the function $C$ must be summable along each geodesic $\pi(u,\varpi)$
in the non-compact case, and have infinite sum along each ray $\pi(u,x)$, where 
$x \in \partial T \equiv \mathcal{X}$ in the compact case, resp.
$x \in \partial^* T \equiv \mathcal{X}$ in the non-compact case.

\smallskip 

The associated hierarchical Laplacian is given for $f$ in the space of test functions
$\mathsf{Loc}_{c}(\mathcal{X})$ of all compactly supported, locally constant
functions on $\mathcal{X}$ by
\begin{equation}\label{eq:lap}
\mathfrak{L}f(x) = \sum_{B \in \mathcal{B}\,:\, B \ni x} C(B)
  \int_B\bigl(f(x)-f(y)\bigr)\, \frac{d\mathsf{m}(y)}{\mathsf{m}(B)}.
\end{equation}
For the following, see e.g. \ocite{Bendikov-Krupski} and 
\ocite{Bendikov-et-al2014},  and note that the 
isotropic Markov semigroups of the latter reference are precisely the ones
induced by hierarchical Laplacians. 

\begin{theorem}\label{thm:hierarchical}
The operator $\mathfrak{L}$ of \eqref{eq:lap}
is densely defined on $L^2(\mathcal{X},\mathsf{m})$ and symmetric, and 
admits a complete system of linearly independent, compactly supported
eigenfunctions. The latter consists of all functions  $f_{B,D}\,$, $D \in \mathcal{B}$ 
with predecessor ball $D'=B$, given by
\begin{equation}\label{eq:fBD}
f_{B,D} = \frac{1}{\mathsf{m}(D)}\boldsymbol{1}_D
- \frac{1}{\mathsf{m}(B)}\boldsymbol{1}_{B}\,.
\end{equation}
The associated eigenvalue is $\lambda(B)$. In addition, when 
$\mathsf{m}(\mathcal{X}) < \infty\,$, one has to add the function $\boldsymbol{1}_X$
with associated eigenvalue $0$. 
\end{theorem}

For distinct $B, \tilde B$, the functions $f_{B,D}$ 
and $f_{\tilde B, \tilde D}$ are orthogonal, while for fixed $B$, 
the finitely many functions $f_{B,D}$ ($D'=B$) are not and can of course 
be orthogonalised. 
In particular, 
$\bigl(\mathfrak{L}, \mathsf{Loc}_{c}(\mathcal{X})\bigr)$ is essentially self-adjoint.

Recall here the tree picture, where balls can be identified with
vertices $u \in T$ , resp. the set of ends ``beyond'' $u$. In this sense, we 
shall also write $f_{B,D} = f_{v,u}$ with $u'=v$,  where $B$ corresponds to $v$ and 
$u$ to $D$. The measure $\mathsf{m}$
is then completely determined by specifying $\mathsf{m}(v)$ such that 
$$
\mathsf{m}(v)=\sum_{u\in T: u'=v} \mathsf{m}(u). 
$$
In the homogeneous case, the natural choice for $\mathsf{m}$ is the unique (up to 
multiplication by constants) measure which is invariant under the isometry group
of $\mathcal{X}$. In terms of the tree, this is the group of all graph automorphisms
which fix each level $S_k\,$, resp. $H_k$ as a set. 
We can normalise by setting $\mathsf{m}(o)=1$  so that 
\begin{equation}\label{eq:Haar}
\mathsf{m}(u) = \mathsf{m}_k = \begin{cases} 
            \dfrac{1}{\mathsf{q}_0\cdots \mathsf{q}_{k-1}}\,,
            &\text{if }\;u \in S_k\,,\;\text{resp.}\; u \in H_k\,,\;k \ge 0\,,\\[9pt]
            \mathsf{q}_{-1} \cdots \mathsf{q}_{-k}\,, &\text{if }\;u \in H_k\,,\;k < 0
            \;\text{ (non-compact case).}
                \end{cases}
\end{equation}
On $\mathcal{X}$, this means of course to specfy $\mathsf{m}(B)$ for all
$B \in \mathcal{B}$.
When $\mathcal{X}$ is itself a totally disconnected locally compact group, 
$\mathsf{m}$ is Haar measure. We write $\mathcal{B}_k$ for the 
set of balls of equal radius $r_k$ corresponding to  the $k^{\text{th}}$ level set 
in the tree,  $B_k(x)$ for the ball in $\mathcal{B}_k$ which contains 
$x \in \mathcal{X}$, and in particular $B_k = B_k(\mathit{0})$, which in the tree
is $\partial T_{o_k}$, resp. $\partial^* T_{o_k}$. 
For a homogeneous hierarchical Laplacian, we also require that the choice function
$C$ is group-invariant, i.e., has constant value $C_k$ on 
each $\mathcal{B}_k\,$. Thus we have a 
\emph{monotone}  sequence of eigenvalues\footnote{In some references, the labelling is
reversed, that is, they have $-k$ in the place of $k$. Here, larger $k$ means 
smaller balls and larger eigenvalues.}
\begin{equation}\label{eq:monotone}
\lambda(B) =\lambda_k \; \text{ for } \; B \in \mathcal{B}_k
\quad \text{with} \quad 0 < \lambda_k < \lambda_{k+1} \quad \text{and}\quad 
\lim_{k \to \infty} \lambda_k = \infty\,.
\end{equation}
Here, $k \in \mathbb{Z}$ in the non-compact case. In the compact case, 
$k \in \mathbb{N}_0$, and there is the additional eigenvalue $\lambda = 0$ with 
eigenfunction $\mathbf{1}_{\mathcal{X}}\,$. The homogenous hierarchical 
Laplacian now reads for $f \in \mathsf{Loc}_{c}(\mathcal{X})$
\begin{equation}\label{eq:homlap}
\mathfrak{L}f(x) =  \sum_{k} C_k \bigl(f(x)-\Pi_kf(x)\bigr)\,, \quad \text{where}
\quad \Pi_kf(x) = \frac{1}{m(B_k)} \int_{B_k(x)} f\, d\mathsf{m}\,.
\end{equation}
Note that having normalised $\mathsf{m}$ as in \eqref{eq:Haar},
$\mathfrak{L}$ is determined by the eigenvalue sequence.

\smallskip

To the authors knowledge, the first operators of the above type to be considered 
were the \emph{spectral multipliers} of \ocite{Taibleson} on $\mathbb{F}^n$, 
where $\mathbb{F}$ is a local 
field. For the sake of simplicity, and since this conveys all the necessary informatin,
we just consider $\mathbb{F} = \mathbb{Q}_{\mathsf{p}}\,$,
where $\mathsf{p}$ is prime. For 
\begin{equation}\label{eq:x}
x = \sum_{k \ge m} a_k\, \mathsf{p}^k \in \mathbb{Q}_{\mathsf{p}}\quad (m \in \mathbb{Z}, 
a_k \in \{0,...,\mathsf{p}-1\}),
\end{equation}
its $\mathsf{p}$-adic norm is $|x|_{\mathsf{p}} = \mathsf{p}^{-m}$ when $a_m \ne 0$, and $|\mathit{0}|_{\mathsf{p}} = 0$. 
The ultrametric on $\mathbb{Q}_{\mathsf{p}}^n$ is induced by the norm (distance to $\mathit{0}$) 
$|(x_1\,,\dots, x_n)|_{\mathsf{p}} = \max \{ |x_j|_{\mathsf{p}} : j=1, \dots, n \}$ for 
$x = (x_1\,,\dots, x_n)$. Thus, each open-compact ball has $\mathsf{p}^n$ successors, 
so that in the associated tree, $\mathsf{q}_k = \mathsf{p}^n$ for all~$k$. 

The Fourier transform on $\mathbb{Q}_{\mathsf{p}}^n$ is induced by any character;
for $x \in \mathbb{Q}_{\mathsf{p}}$ as above we can choose 
\begin{equation}\label{eq:chi}
\chi(x) = \exp{2\pi i \{x\}_{\mathsf{p}}}\,, \quad\text{where} \quad 
\{x\}_{\mathsf{p}} = \sum_{k \le -1} a_k \, \mathsf{p}^k \in \mathbb{Q}
\end{equation}
is the fractional part of $x$. The Abelian group $\mathbb{Q}_{\mathsf{p}}^n$ is 
self-dual, and for $f \in 
L^1(\mathbb{Q}_{\mathsf{p}}^n)$ 
the Fourier transform is
$$
\widehat{f}(\xi) = \int_{\mathbb{Q}_{\mathsf{p}}^r} \chi(x\cdot \xi)\,f(x)\, dx
\,,
$$
where $x\cdot \xi = x_1\xi_1 + \dots + x_n\xi_n$ for $x, \xi \in \mathbb{Q}_{\mathsf{p}}^n$ and 
$dx$ refers to the Haar measure $\mathsf{m}$; see \eqref{eq:Haar}. For the 
inverse transform, we just note that $\;\widehat{\!\!\widehat{f}}(x) = f(-x)$
for $f \in \mathsf{Loc}_{c}(\mathcal{X})$. For $\alpha > 0$, the \emph{Taibleson operator}
$\mathfrak{D}^{\alpha}$ is defined by
$$
\widehat{\mathfrak{D}^{\alpha} f}(\xi) = |\xi|_{\mathsf{p}}^{\alpha} \widehat{f}(\xi)\,,\quad
f \in \mathsf{Loc}_{c}(\mathbb{Q}_{\mathsf{p}}^n).
$$
This operator is a homogeneous hierarchical Laplacian with respect to the Haar measure,
the eigenvalue sequence is $\lambda_k = \mathsf{p}^{k\alpha}$, 
see \ocite{Bendikov-et-al2014}[\S 5.3].  
One computes 
\begin{proposition}\label{pro:TaiVlad} \quad
$\displaystyle \mathfrak{D}^{\alpha} f(x) = \frac{\mathsf{p}^{\alpha}-1}{1-\mathsf{p}^{-\alpha-n}}
\int_{\mathbb{Q}_{\mathsf{p}}^n} \frac{f(x)-f(y)}{|x-y|_{\mathsf{p}}}\, dy.
$
\end{proposition}
For $n=1$, this was introduced as an operator of $\mathsf{p}$-adic fractional derivative
end extensively studied by 
\ocite{Vladimirov-1989}, \ocite{Vladimirov-1994} plus 
many further references, see e.g. the books by
\ocite{Kochubei-2001} and \ocite{Zuniga-2025}.
\\[5pt]
\textbf{B) Spherically symmetric Laplacians.}
In the homogeneous case, the focus is on transition semigroups which are 
invariant under the isometry group of $\mathcal{X}$, which acts doubly transitively 
(transitively on pairs of points at equal distance). 
Let us call them \emph{spherically symmetric,} since they are not exactly the same
as the isotropic processes mentioned above.\footnote{When 
$\mathcal{X} = \mathbb{Q}_{\mathsf{p}}$ is considered, following \ocite{Kochubei-2001}, the 
expression ``rotation invariant'' is used, which makes less sense in more general 
situations like here, whence the author of these notes uses ``spherically symmetric''.}
This means that there is a 
convolution semigroup of probability measures $(\mu_t)_{t > 0}$ on 
$\mathcal{X} \equiv \partial T$ resp. $\equiv \partial^*T$ which  
are determined by the transitions with starting point
$\mathit{0}$. Here $\mathit{0}$ is the neutral elemement of $\mathcal{X}$, 
when the space is identified with an Abelian group; see above, resp. \eqref{eq:group} 
below. 
The $\mu_t$ must have a density with respect to the Haar measure $\mathsf{m}$ of the form 
\begin{equation}\label{eq:iso}
f_t=\sum_k  \frac{c_{k,t}}{\mathsf{m}_k}\mathbf{1}_{B_k}
\end{equation}
with $c_{k,t} >0$, and $p_t(x,\cdot) = \delta_x * f_t\,$. 
Convolution comes from the Abelian group, without depending on the specific
choice of that group acting transitively and freely. 
This is the viewpoint adopted by \ocite{FT-46}.
The other viewpoint, as used by \ocite{FT-41}, is to consider the whole 
isometry group $\mathfrak{G}$ of $\mathcal{X}$, which is locally compact 
and totally disconnected and may be identified with the group of automorphisms
of the tree which stabilise each level set. If $\mathfrak{K}$ is the stabiliser
of $\mathit{0}$, then $\mathcal{X} \cong \mathfrak{G}/\mathfrak{K}$. 
Functions on $\mathcal{X}$ may be identified with functions on $\mathfrak{G}$ which
are constant on the left cosets of $\mathfrak{K}$. Thus, convolution on $\mathfrak{G}$
descends to functions of $\mathcal{X}$. Now in this sense, the density $f_t$ lifted to
$\mathfrak{G}$ is  bi-$\mathfrak{K}$-invariant,
and since $\mathfrak{G}$ acts doubly transitively, the convolution algebra of
bi-$\mathfrak{K}$-invariant functions is commutative: $(\mathfrak{G},\mathfrak{K})$
is a Gelfand pair, see e.g. \ocite{Lang} or more specifically, \ocite{Letac1982}. 
The elements of this convolution algebra are diagonalised by an 
orthogonal  collection of simultaneous eigenfunctions. These are the 
\emph{spherical functions}
of the Gelfand pair. For an introduction and the computation of these functions,
see \ocite{Letac1982}. Relating the notation of \ocite{FT-41} with the one of 
\eqref{eq:fBD}, we have 
$$
\left(\tfrac{1}{\mathsf{m}_k} - \tfrac{1}{\mathsf{m}_{k-1}}\right) \varphi_k = 
f_{B_{k-1},B_k}\,.
$$
In the non-compact case, $k \in \mathbb{Z}$. In the compact case, 
$k\in\mathbb{N}$ and  also the 
constant spherical function $\phi_0 = \mathbf{1}_{B_0}$ with $B_0 = \mathcal{X}$ is 
needed. 
This is the case studied by \ocite{FT-41}. First, a random walk is constructed on 
each level set $S_k$, and then letting $k \to \infty$ and rescaling, a 
$\mathfrak{K}$-invariant convolution semigroup $(\mu_t)_{t > 0}$ on 
$\mathcal{X} \equiv \partial T$ is obtained.  We shall come back to the random walk 
in the next section. For the moment, we only look at the resulting (non-negative
definite) Laplacian $\mathfrak{L}$, where $-\mathfrak{L}$ is the generator of the 
transition semigroup. 

\begin{proposition}\label{pro:compact} {\rm [\ocite{FT-41}]} $\;$
For any choice of the parameter $q >1$, the operator $\mathfrak{L}$ 
satisfies $\mathfrak{L}\varphi_k = q^{k-1}\varphi_k\,$.
Thus, $\mathfrak{L}$ is the homogeneous
hierarchical Laplacian associated with the eigenvalue sequence
$\;\lambda_k=q^k \; \text{for} \; k \ge 0$.\footnote{The values
$q^{k-1}$ given by \ocite{FT-41} come from the fact that the eigenvalue associated with
$\varphi_k$ is $\lambda_{k-1}$ in our numbering.}

In terms of the projections $\Pi_k$ of \eqref{eq:homlap}, 
\begin{equation}\label{eq:FT41}
\mathfrak{L}f = \sum_{k=0}^{\infty} \lambda_k \bigl(\Pi_{k+1}f - \Pi_kf\bigr). 
\end{equation}
\end{proposition}

For $f \in \mathsf{Loc}_{c}(\mathcal{X})$, our space of test functions,
one checks that \eqref{eq:FT41} coincides with the definition of the homogeneous 
hierarchical Laplacian as in \eqref{eq:homlap}. Indeed, $f$ must be a linear combination
of indicator functions of balls, and if $f = \boldsymbol{1}_{B_n(x)}$ then 
$\Pi_{k+1}f - \Pi_kf = f - \Pi_kf = 0$ for $k \ge n$, so that the involved sums range
for $k=0,\dots, n$, and one can rearrange the terms in the sum to
pass from one formula to the other. 

When all forward degrees coincide, this construction covers the Taibleson-Vladimirov
operators on $\mathbb{Z}_{\mathsf{p}}^n\,$, where $\mathbb{Z}_{\mathsf{p}}$ is the ring of 
$\mathsf{p}$-adic integers in $\mathbb{Q}_{\mathsf{p}}$ for arbitrary integer 
$\mathsf{p} \ge 2$.

In the same year, \ocite{Albeverio-1994} published their important 
work on diffusion processes
on $\mathbb{Q}_{\mathsf{p}}\,$, which had a forerunner \ycite{Albeverio-1991}, in turn 
preceded in a slightly different vein 
by \ocite{Evans} .

\smallskip

\textit{One can exclude that at that time, Figà-Talamanca knew the papers by
Albeverio and Karwowski.
There were no internet and no \textsf{arXiv}. Math.Reviews took quite long
after articles had appeared for publishing the review, and the title of 
the proceedings volume in which the forerunner appeared was such that it
is not conceivable that he might have found that paper there. Anyway, it is noteworthy 
that \ocite{Albeverio-1991}, \ycite{Albeverio-1994} cite the older article of
\ocite{FT-1984}. In subsequent 
work, there were mutual citations by Figà-Talamanca and Albeverio,
but the author of these notes does not know whether the two ever met personally.
It also seems that in 1992-93 Figà-Talamanca was not aware of  the paper by
\ocite{Evans}, resp. its not so direct connection to diffusion on compact 
tree boundaries. In any case, the two met at the 1992 conference whose 
proceedings contain the minicourse notes of \ocite{FT-44}, and the later 
work by Del Muto and Figà-Talamanca does refer to Evans.
}

\smallskip

Before returning to a more chronological picture, we anticipate further harmonic
analysis, taken from \ocite{FT-47}, who refer to \ocite{Pontryagin} and \ocite{Evans}. 
If $\mathcal{X}$ (non-compact) 
is identified with an Abelian group $\mathcal{G}$, then there is a double sequence of 
compact-open subgroups $\mathcal{G}_k$ ($k \in \mathbb{Z}$) such that
\begin{equation}\label{eq:group}
\mathcal{G}_k \subset \mathcal{G}_{k-1}\,,\quad \bigcup_k \mathcal{G}_k = \mathcal{G}\,,
\quad \bigcap_k \mathcal{G}_k = \{\mathit{0}\}\,,\;\text{ and }\; 
|\mathcal{G}_k / \mathcal{G}_{k+1}| = \mathsf{q}_k \ge 2.
\end{equation}
A good choice for an associated ultrametric is provided by the Haar measure, normalised
by $\mathsf{m}(\mathcal{G}_0)=1\,$: we set $|\mathit{0}|_{\mathcal{G}}=0$ 
and $|x|_{\mathcal{G}}=\mathsf{m}(\mathcal{G}_k)$ for 
$x \in \mathcal{G}_k \setminus \mathcal{G}_{k+1}$, so that 
$\rho(x,y) = |x-y|_{\mathcal{G}}$.
The cosets of $\mathcal{G}_k$ form 
the level set $\mathcal{B}_k$ of balls, which is $k^{\text{th}}$ horocycle $H_k$ in 
the associated tree, and $\mathsf{q}_k$ is the corresponding forward degree (number
of successors). In particular, $B_k = \mathcal{G}_k\,$.
Let $\Gamma$ be the character group of $\mathcal{G}$. It is the 
union of the compact-open subgroups 
$$
\Gamma_k = \{ \chi \in \Gamma : \chi(x) = 1 \; \text{for all}\; x \in \mathcal{G}_{-k}\}\,,
\quad \text{and} \quad |\Gamma_k/\Gamma_{k+1}| = \mathsf{q}_{-k}\,.
$$
Any Markov semigroup of a standard stationary process on $\mathcal{G}$ is induced
by a convolution semigroup  $(\mu_t)_{t > 0}$ of probability measures as above,
in general not necessarily spherically symmetric, satisfying 
$\lim_{t \to 0} \mu_t = \delta_{\mathit{0}}\,$. Then one has the L\'evy-Khintchine formula, 
see \ocite{Partha}, which in the present setting has a slightly simplified form as in 
\ocite{Kochubei-2001}: there is a positive Radon measure $\mathsf{F}$ on $\mathcal{G}$,
the L\'evy measure, such that 
$\mathsf{a}_k = \mathsf{F}(\mathcal{G} \setminus \mathcal{G}_k) < \infty$ for 
all $k$ and
\begin{equation}\label{eq:levy}
\widehat{\mu}_t(\chi) 
= \exp \left( - t \int_{\mathcal{G}}\bigl(1-\chi(x)\bigr)\, d\,\mathsf{F}(x) \right). 
\end{equation}
In the spherically symmetric case as above in \eqref{eq:iso}, i.e., invariance of $\mu_t$ 
under the full isometry group, the Fourier analysis simplifies further, as proved by
\ocite{FT-46}:

\begin{theorem}\label{thm:DM-FT} {\rm [\ocite{FT-46}]}$\;$
One has the spherical transform 
\begin{equation}\label{eq:spherical}
\langle{\mu}_t\,,\varphi_k \rangle
= \exp \left( - t \int_{\mathcal{G}}\bigl(1-\varphi_k(x)\bigr)\, d\,\mathsf{F}(x) \right). 
\end{equation}
The measure $\mathsf{F}$ has constant density with respect to $\mathsf{m}$ on each set
$\mathcal{G}_k \setminus \mathcal{G}_{k+1} = B_k \setminus B_{k+1}$, so that it is
completely determined by the numbers 
$\mathsf{a}_k = \mathsf{m}(\mathcal{G} \setminus B_k)\,$, which must satisfy 
\begin{equation}\label{eq:ak}
0 \le \mathsf{a}_k \le \mathsf{a}_{k+1} \quad \text{and} 
\quad \lim_{k \to -\infty} \mathsf{a}_k = 0. 
\end{equation}
Any sequence of this type defines a feasible L\'evy measure.
The associated Laplacian satisfies $\mathfrak{L}\varphi_k = \lambda_{k-1}\varphi_k$
with
$$
\lambda_k = \int_{\mathcal{G}}\bigl(1-\varphi_{k+1}(x)\bigr)\, d\,\mathsf{F}(x) 
= \frac{1}{\mathsf{q}_k-1}(\mathsf{q}_k\mathsf{a}_{k+1} - \mathsf{a}_k)\,.
$$
\end{theorem}
For the eigenvalues, see the computations of \ocite{FT-46}[p. 208]. 
We stress that 
Theorem \ref{thm:DM-FT} provides all spherically symmetric Laplacians with
a regular transition semigroup. See also \ocite{Evans}. In the specific case of 
$\mathbb{Q}_{\mathsf{p}}\,$, this had been proved by \ocite{Albeverio-2000}. 

We see that in the spherically symmetric case, we just need the sequence 
$(\mathsf{a}_k)$ satisfying 
\eqref{eq:ak}. In the non-compact case, one should have $\mathsf{a}_k > 0$ for all $k$, and 
in this case, $\lambda_k > 0$, while always $\lim_{k \to -\infty} \lambda_k = 0$,
which allows once more to write 
$$
\mathfrak{L}f = \sum_{k \in \mathbb{Z}} \lambda_k \bigl(\Pi_{k+1}f - \Pi_kf\bigr)\,,
\quad f \in \mathsf{Loc}_{c}(\mathcal{X})
$$
(with $\mathcal{X} = \mathcal{G}$). Note however that the $\lambda_k$ are not necessarily
monotonically increasing. In particular, the spherically symmetric Laplacians and 
associated diffusion 
processes go  beyond the homogeneous hierarchical Laplacians and their isotropic 
Markov processes; compare with \ocite{Bendikov-et-al2014}[Cor. 5.9].

Returning to the historical development, in their important work \ocite{Albeverio-1991},
\ycite{Albeverio-1994}
construct and analyse Markov processes of this type on $\mathbb{Q}_{\mathsf{p}}\,$, and later 
\ocite{FT-46} 
provide the general construction on Abelian groups as above in the first part of their 
paper. We shall come back to the second part in the next section. 
Neither of those references specifies 
the infinitesimal generator directly in the above way. Albeverio and Karwowski use no 
harmonic analysis, but construct their  process starting on horocycles of the tree 
$\equiv$ the collections $\mathcal{B}_k$ of balls, introducing and 
integrating the Kolmogorov equations. Subsequently they analyse the generator as well as
the associated Dirichlet form.\footnote{The author 
of these notes does not agree with the phrase 
in \ocite{Grigoryan2023}[p. 206] which seems to claim that for 
$\mathbb{Q}_{\mathsf{p}}\,$, the class of isotropic Markov processes associated with 
hierarchical Laplacians ``coincides with the ones constructed by Albeverio and 
Karwowski''.}  

Prior to \ocite{FT-46}, \ocite{Albeverio-1999} extended their previous construction on 
$\mathcal{X} = \mathbb{Q}_{\mathsf{p}}$ to a non-group invariant setting. 
Like in the previous material, it is not necessary that $\mathsf{p}$ is prime. 
The tree is such that $\mathsf{q}_k = \mathsf{p}$ on each $H_k$, and there is 
a measure denoted $\rho$ on $\partial^*T \equiv \mathbb{Q}_{\mathsf{p}}$ 
which is  absolutely continuous with respect to the Haar measure given by \eqref{eq:Haar}.
The resulting processes and associated Laplacians rely only on the ultrametric structure 
and no invariance under a group action is assumed. However, besides missing minus signs
after the ``='' in their Kolmogorov equations (2.6ab), they reconsider the distance
transitive case, and their ``Vice versa'' statement 
[Proposition 3.2] appears to be wrong:  it is claimed that monotone eigenvalue 
sequences characterise 
spherically symmetric (``rotation invariant'') Markov semigroups, resp. Laplacians. 
The mistake appears to be the  statement a few line before Prop. 3.2 
that monotonicity of the sequence $\bigl(a(M)\bigr)$  -- here
the L\'evy sequence $(\mathsf{a}_k)$ -- were to imply monotonicity of the eigenvalue
sequence $(h_M)$ (here $\lambda_k$).\footnote{This does by no means affect correctness of 
all the other resuts by Albeverio and Karwowski.} 
In fact, spherically symmetric Laplacians are 
more general than the homogeneous hierarchical ones of \eqref{eq:homlap}. 

Later, \ocite{Albeverio-2008} extended their construction to arbitrary trees
without need of a group action. This was followed by two important papers by
Kigami, see below. In particular, \ocite{Kigami2013} provides one of the most 
general constructions in the spirit of the work described here, but 
this work of Kimgami work has remained somewhat underrepresented in the 
busy activities of  the subsequent 12 years.


\section{Random walks on trees and processes on their boundaries}\label{sec:RW}


The approach of \ocite{FT-41} and \ocite{FT-46}, \ycite{FT-47} as well as of 
\ocite{FT-43}, always in the group-invariant case, consists in starting with
certain nearest neighbour random walks on the tree associated with the ultrametric
space and and using them in a scaling limit which then provides Markovian
semigroups on the boundary of the tree. 

In order to put this into perspective, we first look once more at later work
which relates random walks on trees with boundary processes. This is the significant
work of \ocite{Kigami2010}, \ycite{Kigami2013}, as well as 
\ocite{Bendikov-et-al2014}[\S\S 6--8] which was preceded by
the unpublished note of \ocite{Woess2012}. No group-invariance is assumed there.

A starting point may be the Douglas integral \ycite{Douglas}: the ``energy'' of a
harmonic function on the open unit disk can be written as an integral over the boundary 
(the unit circle) using the Poisson representation of harmonic functions. This transports 
the Dirichlet form of the Brownian motion to a Dirichlet form on the circle -- 
very well explained in the introduction of \ocite{Kigami2010}, where the compact 
case is considered. He then starts with a transient nearest 
neighbour random walk $(Z_n)_{n \ge 0}$ on a tree, where $Z_n$ is the random 
vertex at time $n$. This is a reversible Markov chain with the  
electrical network interpretation, compare with \ocite{WMarkov}[Ch. 4 \& 9]: let 
$P = \bigl( p(u,v) \bigr)_{u,v \in T}$ be the stochastic transition matrix of the walk,
with $p(u,v) > 0$ if and only if $u \sim v$. Then there are positive weights $c(u)$
such that $c(u,v)=c(v,u)$ for all $u,v$ where $c(u,v)=c(u)p(u,v)$. 
The associated Dirichlet form is 
$$
\mathcal{D}_T(g_1\,,g_2) 
= \frac12\sum_{u, v \in T}  c(u,v)\bigl(g_1(u)-g_1(v)\bigr)\bigr(g_2(u)-g_2(v)\bigr)
$$
for functions $g_1\,, g_2: T \to \mathbb{R}$ with finite ``energy'' 
$\mathcal{D}_T(g_i,g_i)$. 
Writing $P^n = \bigl( p^{(n)}(u,v) \bigr)_{u,v \in T}\,$, the Green kernel is
$$
G(u,v) = \sum_{n=0}^{\infty} p^{(n)}(u,v)\,,
$$
so that $0 < G(u,v) < \infty$ for all $u, v$ by transience.
The Martin boundary on which harmonic functions have a Poisson integral representation 
is the boundary of the tree, see \ocite{Cartier}. 
In particular, the random walk $Z_n$ converges almost surely
to a boundary-valued random variable $Z_{\infty}\,$. 
We write $\nu_o$ for the distribution of the latter,
when the starting point is the root vertex $o$. One has $\nu_o(\partial T)=1$ 
and, for $u \ne o$,  
$$
\nu_o(\partial T_u) = F(o,u) \frac{1-F(u,u')}{1-F(u',u)F(u,u')}\,,
$$
where $F(u,v) = G(u,v)/G(v,v)$ is the probability that the walk ever reaches $v$ when its 
starts at $u$. Set $K(u,v) = F(u,v)/F(o,v)$,
the Martin kernel. It has a continuous extension to $\partial T$ given by
$$
K(u,x) = K(u, u \wedge x) \,, \quad K(o,x)=1.
$$
Then every bounded harmonic function $h$  has a unique integral representation 
$$
h(u) = h_{f}(u) = \int_{\partial T} f(u) K(u,x) \, d\nu_o(x)\,,
$$
where $f \in L^{\infty}(\partial T,\nu_o)$. An analogous represenation 
(with $f$ not necessarily bounded) also holds for every harmonic function 
with finite energy. Let us assume here that $F(v,v') < 1$ for all $v \in T$, so that 
$\nu_o$ is supported by the whole of $\partial T$.  We may also assume that 
the Green kernel 
vanishes at infinity, so that $h_{f}$ provides the continuous extension on 
$T \cup \partial T$ of every continuous function on $\partial T$.\footnote{This second 
assumption is not used by Kigami, but anyway, $\lim_{v\to x} G(v,o) = 0$ for 
$\nu_o$-almost
every $x \in \partial T$, see \ocite{WMarkov}[Thm. 9.43].} 
In view of the above, $\mathcal{D}_T$ gives rise to the Dirichlet form
$$
\mathcal{D}_{\partial T}(f_1\,,f_2)  = \mathcal{D}_{T}(h_{f_1}\,,h_{f_2}).
$$
\ocite{Kigami2010} computes its integral kernel and shows that the form is
regular on $L^2(\partial T, \nu_o)$: by the theory of Dirichlet forms
-- see \ocite{Fukushima} -- it induces a symmetric diffusion process on $\partial T$.
As a matter of fact, that integral kernel is the \emph{Na{\"\i}m kernel} of
\ocite{Doob}, adapted to the setting of reversible Markov chains by 
\ocite{Silverstein}. In the case of random walks on trees, it has a simple 
expression \& proof, see \ocite{Bendikov-et-al2014}[\S 6]:

\begin{proposition}\label{pro:naim} The Na{\"\i}m kernel on $\partial T$ is
$$
\Theta_o(x,y) = \begin{cases} \dfrac{c(o)}{G(o,o)F(o,x \wedge y)F(x \wedge y,o)} &
                              \text{if } x \ne y\,,\\[5pt]
                              +\infty & \text{if } x = y\,.
                \end{cases}
$$
With the reference measure $\nu_o\,$, the associated Laplacian is given on 
$\mathsf{Loc}_{c}(\mathcal{\partial T})$ by 
$$
\mathfrak{L}f(x) = \int_{\partial T} \bigl(f(x)-f(y)\bigr)\Theta_o(x,y)\, d\nu_o(y).
$$
It is the hierarchical  Laplacian on $L^2(\partial T, \nu_o)$ with eigenfunctions
$f_{B,D} = f_{v,u}$ with $u'=v$ as in \eqref{eq:fBD}, where $B = \partial T_v$ 
and $D = \partial T_u\,$. The eigenvalues are
$\lambda(v) = 1/G(v,o)$, $v \in T$. In addition, there is the constant function 
$\mathbf{1}$ with eigenvalue $0$.
\end{proposition}

\ocite{Kigami2013} then studies the analogous construction in the non-compact case,
where $\mathcal{X}=\partial^* T$  by moving the reference point to $\varpi$. 
The above can be done for any ``root'' vertex in the place of $o$, in
particular for the $n^{\text{th}}$ predecessor $o_n$ of $o$. Define the Radon measure
$\mathsf{m}$ on $\partial^*T$ for compact $A \subset \partial^* T$ and the rescaled 
Na{\"\i}m kernel $J$ for $x,y \in \partial^* T$ by 
$$
\mathsf{m}(A) = \lim_{n \to \infty}
\frac{\nu_{o_n}(A)}{\nu_{o_n}(\partial^* T_o)}
\quad \text{and} \quad 
\lim_{n \to \infty} J(x,y) = \Theta_{o_n}(x,y) \bigl(\nu_{o_n}(\partial^* T_o))\bigr)^2
$$
Both sequences stabilise for $n \ge k$ when $A \subset \partial^*T_{o_k}\,$, resp.
$x \curlywedge y \in \partial^*T_{o_k}\,$, see \ocite{Bendikov-et-al2014}[\S 7.3]. 
One can then write $J(x,y) = \mathfrak{j}(x \curlywedge y)$ with
$\mathfrak{j}(v) = \vartheta^2 G(v,v)\big/\bigl(K(v,\varpi)^2 c(v)\bigr)$, 
where $\vartheta > 0$ is known. If the Green kernel vanishes at infinity then 
for $f_1\,,f_2 \in \mathsf{Loc}_{c}(\partial^*T)$, the Dirichlet form
$\mathcal{D}_{\partial T}(f_1\,,f_2)$ can be rewritten with integral kernel
$J(\cdot,\cdot)$ with respect to the measure $\mathsf{m}$. The associated
Laplacian is
$$
\mathfrak{L}f(x) = \int_{\partial T} \bigl(f(x)-f(y)\bigr)J(x,y)\, d\mathsf{m}(y)\,,
$$
the eigenfunctions are again all $f_{v,u}$ with $u'=v$ and associated
eigenvalue $\lambda(v) = K(v,\varpi)/\vartheta\,$, $v \in T$.

In \ocite{Kigami2010}, \ycite{Kigami2013} as well as \ocite{Bendikov-et-al2014}, 
it is shown
that conversely, every hierarchical Laplacians on $\partial T$, resp. $\partial^*T$
and the associated transition semigroup are induced in the above way by a random walk,
up to a possible time change by a constant factor.

In the group-invariant case, $\nu_o$ and $\mathsf{m}$ are the respective Haar measures,
the jump kernels $\Theta_o$ and $J$ are group-invariant, and the eigenvalues $\lambda(v)$
depend only on the level set index of $v$, so that one gets all the 
homogeneous hierarchical Laplacians and associated diffusion processes. 

\smallskip

Now we finally come to the core of the constructions of \ocite{FT-41}, \ocite{FT-43} and 
\ocite{FT-46}, \ycite{FT-47}. There, the respective
constructions also rely on discrete time random walks on the tree asociated with the
ultrametric space, resp. totally disconnected Abelian group, but in a different way 
via the convolution semigroup $(\mu_t)_{t > 0}\,$.
\\[5pt]
\textbf{A) \ocite{FT-41}.} 
In the compact, homogeneous case, the author starts with a probability parameter
$p \in (0\,,\,\frac12)$ (the degenerate case $p = \frac12$ is also considered
but less significant; $p \in (\frac12\,,\,1)$ does not lead to a diffusion).
The random walk moves from the root $o$ to any of its neighbours $u$ with probability
$p(o,u)=1/\mathsf{q}_0\,$. For a vertex $u \in S_k$ ($k \ge 1$), the random walk moves up 
to $u'$ with probability $p(u,u') = p$ and down to any of the successors $v$ of $u$ with 
probability $p(u,v)= (1-p)/\mathsf{q_k}\,$. This walk is transient.

Now the tree is truncated at level $S_n\,$, removing the vertices $v$ with $|v| > n$. 
For $u \in S_n$ the transition probabilities $p(u,v)$ are adapted such that the walk 
moves to  $u'$ with probability $1$. 
Then a Markov chain on $S_n$ is constructed: for $u, v \in S_n\,$, $p_n(u,v)$
is the probability that the random walk on the truncated tree, starting at $u$, first 
re-enters $S_n$ at $v$. This is computed as follows: set $q = (1-p)/p$ and recall
$\mathsf{m}_k$ from \eqref{eq:Haar}.
$$
p_n(u,v) = 
\sum_{j= \max\{1,k\}}^{n-1} \left(\frac{q-1}{q^j-1} - \frac{q-1}{q^{j+1}}\right) 
\frac{\mathsf{m}_n}{\mathsf{m}_{n-j}} + \frac{q-1}{q^n-1}\mathsf{m}_n\,,
\quad \text{if }\; |u \wedge v| = n-k\,.
$$
We choose the reference point $o_n \in S_n$ and define the probability measure
$\sigma_n$ on $\mathcal{X}$ with density\footnote{It is denoted $\mu_n$ by
Figà-Talmanca; in order to avoid confusion with $\mu_t$ we write $\sigma_n$ here.}
\begin{equation}\label{eq:fn}
f_n(x) = \sum_{u \in S_n} \frac{p_n(o_n,u)}{\mathsf{m}_n} \mathbf{1}_{\partial T_u}  
\end{equation}
with respect to the measure of \eqref{eq:Haar}. 
Thus, in the Gelfand pair approach outlined above, or equivalently if $\mathcal{X}
= \mathcal{G}$ is a totally disconnected Abelian group, $\sigma_n$ is a 
spherically symmetric probability measure, and 
$p_n(go_n\,, ho_n) = \sigma_n(g^{-1}h\,\partial T_{o_n})$.
Its eigenvalues as a convolution operator are computed via the spherical functions:
$$
\langle \sigma_n\,,\varphi_k \rangle = \begin{cases} 1\,,& k=0\,,\\
                                    1- \dfrac{q-1}{q^{n-k+1}-1}\,, & k=1, \dots, n-1\,,\\  
                                    0\,, & k \ge n\,. 
                                    \end{cases}
$$
The ones up to $k=n$ are the eigenvalues of $\bigl(p_n(u,v)\bigr)_{u,v \in S_n}$.

\begin{theorem}\label{thm:FT41} {\rm [\ocite{FT-41}]} $\;$
Let $k(n) = (q^n-1)/(q-1)$, so that for the second largest eigenvalue of 
$\bigl(p_n(u,v)\bigr)_{u,v \in S_n}\,$, 
$$
\langle \sigma_n\,,\varphi_1 \rangle^{\lfloor tk(n)\rfloor} \to e^{-t} \quad 
\text{from below,}
$$
where $\lfloor \cdot \rfloor$ stands for the next lower integer.
Then for any $t > 0$, the sequence of convolution
powers
$$
\sigma_n^{(\lfloor tk(n)\rfloor)} \to \mu_t\,,\quad \text{as }\; n \to \infty\,,
$$
where $(\mu_t)_{t > 0}$ is the convolution semigroup which leads to the 
Laplacian of \eqref{eq:FT41} with $\lambda_k = q^k$.
\end{theorem} 

We compare this with the boundary process obtained via the Na{\"\i}m kernel
from the initial random walk with transition matrix $\bigl(p(u,v)\bigr)$ described above:
in our spherically symmetric situation, $\nu_o = \mathsf{m}$ is the same measure as 
here (equidistribution on $\partial T$ as seen from the root $o$), the eigenfunctions are 
the same, and for $v \in S_k$, one easily computes $G(v,o) = q^{1-k}/(q-1)$,
whence the eigenvalue sequence in our numbering is $\lambda_k = (q-1) q^k$, 
Thus, the two processes are the same up to the time change by the constant
factor $q-1\,$: if $(\mu_t)$ is the convolution semigroup constructed by 
\ocite{FT-41} and $\tilde \mu_t$ is the one corresponding to the boundary process, then 
$\tilde\mu_t = \mu_{(q-1)t}\,$. 

\textit{Time ago, A. Figà-Talamanca pointed out a few misprints to the author of 
these notes: 
in \ocite{FT-41}, on page 162, lines 6$-$ and 5$-$
and page 164, lines 8 and 9, it should be $p \ne 1/2$ and $p = 1/2$ (instead of $2$). On 
page 166, L3, the last term should be divided by $\rho_n-1$ in the notation of that 
paper.} 
\\[5pt]
\textbf{B) \ocite{FT-43}.} 
The three authors apply the same method as above to the non-compact case  in the 
first part of this work. Here $\mathcal{X}$ is a local field, and the associated tree
in horocyclic (``hierarchical'') ordering is such that $\mathsf{q}_k = \mathsf{q}$
for all $k \in \mathbb{Z}$. For the sake of simplicity, let us consider 
$\mathbb{Q}_{\mathsf{p}}\,$, so that $\mathsf{q} = \mathsf{p}$, but this is 
not relevant, $\mathsf{q}$ can be any integer
$\ge 2$. The aim is to 
construct $\alpha$-stable random variables on $\mathbb{Q}_{\mathsf{p}} = \partial^* T\,$,
where  $\alpha > 0$. Call the random variable $U$, write $\mu$ 
for its distribution and recall the character $\chi$ from \eqref{eq:chi} in the 
case of dimension $n=1$ (arbitrary $n$ works equally well). 
Then $U$, resp. $\mu$ are required to satisfy
\begin{equation}\label{eq:stable}
\mathbf{E}\bigl(\chi(\xi U)\bigr) = \widehat{\mu}(\xi) 
= e^{-c|\xi|^{\alpha}} \quad \text{for all }\; \xi \in 
\widehat{\mathbb{Q}}_{\mathsf{p}} = \mathbb{Q}_{\mathsf{p}}
\end{equation}
for some constant $c > 0$ whose value is not relevant. Here, $\mathbf{E}$ denotes 
expectation. It is noteworthy that there is a non-trivial solution for any 
$\alpha > 0$, while in the case of $\mathbb{R}$ it must be $\alpha \le 2$.
The construction is analogous to the above compact case. Namely, $\mu = \mu_1$
in a convolution semigroup $(\mu_t)_{t > 0}$ on $\partial^* T$. The authors start with 
the nearest neighbour random walk on $T$ where for predecessors with respect to $\varpi$,
$$
p(u,u')  = \frac{1}{1+\mathsf{p}^{\alpha}} \quad \text{and} \quad 
p(u',u)  = \frac{\mathsf{p}^{\alpha}}{\mathsf{p}(1+\mathsf{p}^{\alpha})}.
$$
The tree is truncated at $H_n\,$, removing the vertices $v$ with $\mathsf{hor}(v) > n$.
The outgoing transition probabilities $p(u,v)$ at $u \in H_n$ are adapted  
such that the walk moves to  $u'$ with probability $1$. Then the induced walk on $H_n$ 
is considered, for which $p_n(u,v)$ is the probability that starting at $u$, the first
return to $H_n$ occurs at $v$. It is computed similarly to the compact case, and 
$p_n(o_n,\cdot)$ is lifted to a spherically symmetric probability measure $\sigma_n$ 
with density $f_n$ on 
$\mathbb{Q}_{\mathsf{p}}$ in the same way as in \eqref{eq:fn} (with $H_n$ and 
$\partial^* T_u$ in the place of $S_n$ and $\partial T_u\,$, respectively).
Then the authors compute
$$
\langle \sigma_n\,,\varphi_k \rangle = \langle f_n\,,\varphi_k \rangle
= \begin{cases} 
1- \dfrac{\mathsf{p}^{\alpha}-1}{\mathsf{p}^{(n+1-k)\alpha}-1}\,, & k<n\,,\\  
0\,, & k \ge n\,. 
\end{cases}
$$
Corresponding to the $k(n)$ from above, here 
$j(n) = \mathsf{p}^{(n+1)\alpha}/(\mathsf{p}^{\alpha}-1)$, and 
the sequence of convolution powers $\sigma_n^{(\lfloor j(n) \rfloor)}$ 
converges to the desired  measure $\mu$. This is the substance; 
the presentation of the authors
is more  probabilistic: 

\begin{theorem}\label{thm:FT43}{\rm [\ocite{FT-43}]} $\;$
Let $(V_k)_{k \in \mathbb{N}}$ be a sequence of i.i.d. random 
variables on $\mathbb{Q}_{\mathsf{p}}$ with distribution $\sigma_0\,$. Then 
$\mathsf{p}^n V_k$ 
(product in $\mathbb{Q}_{\mathsf{p}}$) has distribution $\sigma_n$, and 
$$
U_n = \mathsf{p}^n \sum_{k=1}^{j(n)} V_k
$$
converges in distribution to a random variable. Its law $\mu$ has density 
$$
f(x) = \frac{\mathsf{p}-1}{\mathsf{p}}\sum_{k=-\infty}^{\infty}  
\mathsf{p}^k e^{-\mathsf{p}^{k\alpha}}\varphi_k(x)\,,
$$
with respect to the Haar measure, and $\widehat\mu$ satisfies \eqref{eq:stable}.
\end{theorem}

If instead one considers 
$\mu_t = \lim_n \mu_n^{(\lfloor tj(n) \rfloor)}$, then it is a diffusion
semigroup, and at this point it will be no surprise that this is the one associated with 
the Taibleson-Vladimirov operator of Proposition \ref{pro:TaiVlad}. 
Once more, the boundary process obtained via the rescaled Na{\"\i}m kernel $J(x,y)$
from the original random walk with transition matrix $\bigl(p(u,v)\bigr)$ on the 
tree associated with $\mathbb{Q}_{\mathsf{p}}$ is the same process as the one 
with semigroup $(\mu_t)$ up to a time change by a constant factor, compare with 
\ocite{Bendikov-et-al2014}[\S 8.3].

As the authors state in the introduction (and the reader should understand
after all that has been explained so far), the same methods apply to arbitrary 
local fields $\mathbb{F}$, the vector spaces  
$\mathbb{F}^n$, as well as -- with small modifications -- to general totally 
disconnected Abelian groups as in \eqref{eq:group}. A further generalisation is 
mentioned at the end of \ocite{FT-47}.

\smallskip

\textit{We remark that \ocite{FT-43} cite \ocite{Taibleson}, but not 
\ocite{Vladimirov-1989}. Indeed, Mitchell H. Taibleson was a well known harmonic analyst
and frequent guest in Rome, while apparently Figà-Talamanca was not aware of the 
work of Vladimirov. In turn, \ocite{Bendikov-et-al2014} do not cite any of
the present work of Figà-Talamanca, an omission (overcome in later work) for which 
the author of these notes assumes responsibility.
}

\smallskip

Let us also mention the second, shorter part of the paper by \ocite{FT-43}: 
it is not pertinent
to ultrametric analysis, but it uses the far-reaching analogy between horocyclically 
ordered (in particular, regular) trees and the hyperbolic upper half plane. 
The authors show how hyperbolic Brownian motion with downward drift lead to 
limit distributions on the real line (here the geometric analogue of $\partial^*T$) 
which are  in the domain of attraction of $\alpha$-stable random variables.
\\[5pt]
\textbf{C) \ocite{FT-46}.} 
We now return to the most cited one of the papers of Figà-Talamanca commented upon here.
For their construction in the group-invariant, non-compact case, the authors also use a 
nearest neighbour random walk on the tree, but its transition probabilities are not $>0$
along all edges. Recall that we denote by $o_k$ ($k \in \mathbb{Z}$) the vertex in $H_k$ 
on the bi-infinite geodesic $\pi(\varpi,\mathit{0})$  from $\varpi$ to 
$\mathit{0} \in \partial^*T$. Let $T_k$ be the subtree of $T$ which is spanned by its 
root at $o_k$ and all vertices which have $o_k$ as an ancestor (iterated
predecessor), but not 
$o_{k+1}\,$. The $T_k$ are disjoint, and their union is all of $T$, comprising all
edges except those which lie on the geodesic $\pi(\varpi,\mathit{0})$. 
Then a random walk is 
constructed as follows. We may start with the sequence $(\mathsf{a}_k)_{k \in \mathbb{Z}}$
of \eqref{eq:ak} which defines a spherically symmetric L\'evy measure, and here we
assume that $\mathsf{a}_k > 0$ strictly and $\mathsf{a}_k \to \infty$ as 
$k \to \infty\,$. Then we let $\mathsf{f}_k = \mathsf{a}_k/\mathsf{a}_{k+1}$ 
(written as $\zeta_k\,$, resp. $\zeta_i$ by Del Muto and 
Figà-Talamanca). The numbers are used to define a nearest neighbour random walk for which 
\begin{equation}\label{eq:rwtk}
\begin{aligned}
&p(u, u') = \frac{\mathsf{f}_k}{1+\mathsf{f}_k} \; \text{ for }\; u \in T_k\,,\quad 
p(u',u) = \frac{1}{\mathsf{q}_{\mathsf{hor}(u')}}\!\cdot\!\frac{1}{1+\mathsf{f}_k}  \; 
             \text{ for }\; u \in T_k\,,\; u' \ne o_k\,, \\
&p(o_k\,,o_{k+1})=0\,,\hspace*{2.8cm}
p(o_k,u) = \frac{1}{\mathsf{q}_k-1}\!\cdot\!\frac{1}{1+\mathsf{f}_k}  \; 
             \text{ for }\; u \in T_k\,,\; u' = o_k\,.
\end{aligned}
\end{equation}
The numbers are chosen such that for any $u \in T_k$ one has 
$\mathsf{f}_k = F(u,u')$, the probability that the walk starting at $u$ ever reaches the
predecessor $u'$. We now fix $n$ and construct a probability measure 
$\mu_n$ on $\mathcal{X}= \partial^*T$.
\ocite{FT-46} describe this via a discrete probability measure on $H_n\,$. 
Contrary to the preceding constructions, at $u \in H_n \setminus \{ o_n\}$ 
one considers the probability that for the above random walk on $T$ 
starting at $o_{n-1}\,$, the {\emph{last} (instead of first, as before) visit to  $H_n$ 
occurs at $u$. (If we want the random walk
to start at $o_n$, then we may truncate the tree $H_n$ as previously and redefine 
the transition probability from $o_n$ to $o_{n-1}$ to be $=1$.) 

Here we use an explanation slightly different from the authors'.
First observe that by the assumption that $\mathsf{a}_k \to \infty$
there are arbitrarily large $k$ such that $\mathsf{f}_k < 1$. 
Since $p(o_k\,,o_{k+1}) = 0$,
the walk is transient, so that for any starting point, $Z_n$ 
will converge almost surely to
a $\partial T$-valued random variable $Z_{\infty}\,$. 
Its distribution, when the starting point
is $o_{n-1}\,$, is denoted $\nu_n$. Then the probability that 
$Z_{\infty}$ is in $\partial^* T_{o_k} \setminus \partial^* T_{o_{k+1}}$ (the set 
of ends of $T_k$) is\footnote{By \ocite{FT-46}, 
this is written as $\lambda_n(U_{n-m})$ with $n-m=k$.}
\begin{equation}\label{eq:nun}
\nu_n(\partial^* T_{o_k} \setminus \partial^* T_{o_{k+1}}) 
= \frac{\mathsf{a}_{k+1}-\mathsf{a}_k}{\mathsf{a}_n}\,,
\quad k \le n-1,
\end{equation}
because the probability to reach $o_{k-1}$ from $o_n$ is 
$\mathsf{f}_{n-1} \mathsf{f}_{n-2} \cdots \mathsf{f}_{k} 
= \mathsf{a}_k/\mathsf{a}_n\,$,
and thereafter the walk cannot return to $o_{k-1}$.
In particular, we see that 
$$
\nu_n(\partial^* T)=
\sum_{k \le n-1}\nu_{o_{n-1}}(\partial^* T_{o_k} \setminus \partial^* T_{o_{k+1}}) =1\,,
$$
so that the limit random variable is indeed $\partial^*T\setminus \{\mathit{0} \}$-valued. 
Then the above-mentioned ``last exit'' probability distribution on 
$H_n \setminus \{o_n\}$
is lifted to the spherically symmetric probability measure 
$\sigma_n$ on $\partial^*T$ with density 
$$
f_n = \sum_{u \in H_n \setminus \{o_n\}} \frac{\nu_n(\partial^* T_u)}{\mathsf{m}_n} 
\, \mathbf{1}_{\partial^* T_u}
$$
analogous to \eqref{eq:fn}.
($\sigma_n$ is denoted $\lambda_n$ by the authors.)
Note that $\sigma_n(\partial^* T_u)$ is the same for each of the 
$(\mathsf{q}_k-1)\mathsf{q}_{k+1}\cdots \mathsf{q}_{n-1} 
=(\mathsf{m}_k - \mathsf{m}_{k+1})/\mathsf{m}_n$ 
elements $ u \in T_k \cap H_n\,$,  where $\mathsf{m}_k$ is as in \eqref{eq:Haar}.
We can now see $\sigma_n$ on $\mathcal{X}$  as a convolution operator in the Gelfand pair 
setting, or the equivalent setting where $\mathcal{X}=\mathcal{G}$ is a totally 
disconnected Abelian group as in \eqref{eq:group}. 
The spherical transform of $\sigma_n$ is computed as
$$
\langle \sigma_n\,,\varphi_k \rangle 
= \begin{cases} 
1- \dfrac{1}{\mathsf{a}_n} \cdot \dfrac{\mathsf{q}_{k-1}\mathsf{a}_k 
                      - \mathsf{a}_{k-1}}{\mathsf{q}_{k-1}-1}  \,, & k<n\,,\\  
0\,, & k \ge n\,. 
\end{cases}
$$
\begin{theorem}\label{thm:FT46} {\rm [\ocite{FT-46}]} $\;$
For each $t > 0$, the sequence of convolution powers 
$\sigma_n^{(\lfloor t\mathsf{a}_n \rfloor)}$ converges to a spherically symmetric
probability measure $\mu_t$ on $\mathcal{X}\,$, and $(\mu_t)_{t > 0}$ is 
the convolution semigroup with satisfies \eqref{eq:spherical} for the spherically
symmetric L\'evy measure defined by the sequence $(\mathsf{a}_k)_{k \in \mathbb{Z}}\,$.
\end{theorem}
This construction leads to all spherically symmetric diffusion semigroups in the 
homogeneous case for which the L\'evy sequence is non-degenerate in the sense 
that $\mathsf{a}_k >0$ for all $k$ and $\lim_{k \to \infty} \mathsf{a}_k = \infty\,$,
including the isotropic processes associated with hierarchical Laplacians. 
\\[5pt]
\textbf{D) \ocite{FT-47}.}
This last paper is short, just displaying three lemmas besides the explanatory 
text, but it also provides noteworthy simplifications of the author's previous work.  
In this case, the ultrametric space is
a totally disconnected Abelian group $\mathcal{G}$ as in \eqref{eq:group}. The diffusion 
processes under consideration are $\mathcal{G}$-invariant, whence once more 
given by a convolution semigroup $(\mu_t)_{t > 0}$ with $\lim_{t \to 0} \mu_t = 
\delta_{\mathit{0}}\,$. But now the probability measures $\mu_t$ are not
assumed to be spherically symmetric (``non-isotropic'' in the authors' terminology). 
Thus, the L\'evy representation \eqref{eq:levy} does not simplify to the spherical 
transform. Nevertheless, the authors show that there is an approximation 
induced by a random walk which generalises \eqref{eq:rwtk}.
It is again assumed that the L\'evy measure $\mathsf{F}$ satisfies $\mathsf{a}_k > 0$ and
$\lim_{k \to \infty} \mathsf{a}_k = \infty\,$, where 
$\mathsf{a}_k = \mathsf{F}(\mathcal{G} \setminus \mathcal{G}_k)$. The subtrees $T_k$ are as
above. If $u \in H_k \setminus \{ o_k \}$, then the ball 
$\partial^* T_u \subset\partial^*T=\mathcal{G}$ is of the form $x + \mathcal{G}_k\,$,
$x \in \mathcal{G} \setminus \mathcal{G}_k\,$, and we write $\mathsf{F}(u)$
for $\mathsf{F}(x + \mathcal{G}_k)$. Also, by a small abuse of this notation, 
let us write 
$\mathsf{F}(o_k) = \mathsf{F}(\mathcal{G}_k \setminus \mathcal{G}_{k+1}) 
= \sum_{u \in T_k : u'=o_k} \mathsf{F}(u) = \mathsf{a}_{k+1} - \mathsf{a}_k\,$.
With $\mathsf{f}_k = \mathsf{a}_k/\mathsf{a}_{k+1}$ as defined above, the probabilities
of the random walk on $T$ now are 
\begin{equation}\label{eq:rwtkF}
p(u, u') = \frac{\mathsf{f}_k}{1+\mathsf{f}_k}\,,\quad
p(u',u) = \frac{\mathsf{F}(u)}{\mathsf{F}(u')}\cdot\frac{1}{1+\mathsf{f}_k}  
\; \text{ for }\; u \in T_k\,, \quad p(o_k\,,o_{k+1})=0\,.
\end{equation}
The definition has to be adapted when $F(u') = 0$. In this case,
the probability $p(u',u)$ is defined as in \eqref{eq:rwtk}. Then,
as before, the random walk converges 
to a $\partial^*T \setminus \{\mathit{0} \}$-valued random 
variable, and we let $\nu_n$ be its distribution when the starting point is
$o_{n-1}\,$. Again, it is written as $\lambda_n$ by the authors who then show the 
following.
\begin{lemma}\label{lem:FT47} For $u \in H_n \setminus \{o_n\}\;$ one has 
$\; \nu_n(\partial^* T_u) = \mathsf{F}(u)/\mathsf{a_n}\,.$
\end{lemma}
{\it Remark: the proof of this in \ocite{FT-47}[p. 227] is a bit incomplete. Therefore 
a more complete outline is added here. 
\\[3pt]
{\rm Proof.} There is $k<n$ such that $u \in T_k\,$.
\\[3pt]
{\rm (1)} If 
$\mathsf{F}(o_k) = \mathsf{a}_{k+1} - \mathsf{a}_k = 0$ then also $\mathsf{F}(u) =0$.
Also, $\nu_n(\partial^* T_u)  =0$ by \eqref{eq:nun}. The satement of the lemma 
holds in this case.
\\[3pt]
{\rm (2)} If $\/\mathsf{F}(o_k) > 0$ and $u \in S_n \cap T_k$ is such that 
$\mathsf{F}(u) > 0$ then also $\mathsf{F}(v) > 0$ for all 
$v \in \pi[o_k\,,u] \setminus \{ o_k \}$. For such $v$, one has 
the intuitively clear, but not completely obvious identity
$\nu_n(\partial^* T_{v})/\nu_n(\partial^* T_{v'}) =\mathsf{F}(v)/\mathsf{F}(v')$.
This follows from \eqref{eq:rwtkF}, the fact that $F(v,v') = \mathsf{f_k}$ and 
the formula given by \ocite{WMarkov}[Cor. 9.57]. Together with the 
three line computation  of \ocite{FT-47}[p. 227] this shows that stated identity holds. 
\\[3pt]
{\rm (3)} If $\/\mathsf{F}(o_k) > 0$ and $u \in S_n \cap T_k$ is such that 
$\mathsf{F}(u) = 0$ then there is $v \in \pi[o_k\,,u] \setminus \{ o_k \}$
such that $F(v) = 0$ but $F(v') \ne 0$. Then $v \in H_m$ for some $m \le n$.
For each element of the non-empty set 
$A = \{ \tilde u \in H_m : \tilde u' = v'\,,\; \mathsf{F}(\tilde u) > 0\}$ 
we have by the above
$\nu_n(\partial^* T_{\tilde u})/\nu_n(\partial^* T_{v'}) 
=\mathsf{F}(\tilde u)/\mathsf{F}(v')$.
Therefore 
$$
\nu_n(\partial^* T_{v'}) = \sum_{\tilde u\in A}\nu_n(\partial^* T_{\tilde u}).
$$
But this implies $\nu_n(\partial^* T_{v}) = 0$, whence also $\nu_n(\partial^* T_{u}) = 0$,
and the formula holds also in ths case.}
%

\smallskip

At this point, a surprisingly simple argument of Harmonic Analysis of
\ocite{FT-47}[Lemma 2] shows that the formula of Lemma \ref{lem:FT47}
implies the following.

\begin{theorem}\label{thm:FT47} For the sequence of probability distributions $\nu_n$
associated with the random walk \eqref{eq:rwtkF}, one has that for every $t >0$, the 
convolution powers $\;\nu_n^{((\lfloor t\mathsf{a}_n \rfloor)}$ converges to  $\mu_t\,$,
where $(\mu_t)_{t > 0}$ is the convolution semigroup on the totally disconnected
Abelian group $\mathcal{G}$ associated with the L\'evy measure $\mathsf{F}$.
\end{theorem}

In the final section, it is explained under which conditions the process is
stable, generalising \ocite{FT-43}. 

Note that in \ocite{FT-46}, \ycite{FT-47}, the random walks are quite different
from those of the preceding papers.  They rely strongly on the homogeneity of the 
ultrametric space and there is no apparent duality bewtween random walks
and boundary processes as via the Na{\"\i}m kernel in the spirit of the 
Douglas integral. Theorem \ref{thm:FT47} both generalises and simplifies 
Theorem \ref{thm:FT46} from the preceding paper. It shows how arbitrary diffusion
process with unbounded L\'evy measure on those groups can be obtained from 
nearest neighbour random walks on the underlying tree, and also the specification
of the L\'evy measure is associated with the tree via the values $\mathsf{F}(u)$.

\smallskip 

In conclusion, the principal interest of Figà-Talamanca in all 4 papers
commented upon here is to work out how the combination of harmonic analysis 
and different types of discrete-time nearest neighbour random walks on trees 
can be used to
construct, resp. approximate and understand diffusion processes on 
homogeneous ultrametric spaces. The focus is not on the subsequent analysis 
such as heat kernel asymptotics, Schr\"odinger equations and other issues in 
the spirit of partial differential equations. In particular, in this respect the 
bibliography of the present notes is not intended to be complete. 

\textit{Regarding the citations received by the present work of Figà-Talamanca,
on one hand it may be true that he himself sometimes started
to study the existing literature in this area only after having
already elaborated his own approach. On the other hand, it appears that 
his contributions, although possibly out of the main schools, must have 
been considered well related to them and appreciable, given the fact 
that, in their subsequent bodies of work, leading mathematicians like 
Albeverio and Karwowski as well
as Bendikov and Grigor{\cprime}yan and  their respective
collaborators have frequently included the work of Figà-Talamanca in 
their citation lists. The way each mathematician selects his quotes is 
perhaps based on quite peculiar decisions and choices, and it may be 
difficult to acknowledge contributions from authors outside the 
established schools.
Indeed, the author of these notes found a bit 
surprising that -- with a single exception -- the articles on ultrametric 
issues by Figà-Talamanca and collaborators were never quoted
in the work of the author of all their Mathscinet reviews.} 

\section*{Relevant bibliographic sources}\label{sec:bibliography}

\newpage


\section*{Citations of Figà-Talamanca's papers 
commented in this chapter}\label{sec:citations}

\noindent
{\small 
FT-1994 stands for~\ocite{FT-41}*{Proceedings,~San~Antonio~1993},\\ 
B-CT-FT-2001 stands for~\ocite{FT-43}*{Pacific~J.~Math.}\\
DM-FT-2004 stands for~\ocite{FT-46}*{Expos.~Math.},\\ 
DM-FT-2006 stands for~\ocite{FT-47}*{Pacific~J.~Math.}\\

\end{document}